\documentclass[11pt, final]{article}
\usepackage{a4}
\usepackage{amsmath}%
\usepackage{amstext}%
\usepackage{amssymb}%
\usepackage{showkeys}%
\usepackage{epsfig}%
\usepackage{cite}
\newtheorem{theorem}{Theorem}

\newtheorem{example}[theorem]{Example}

\newtheorem{proposition}[theorem]{Proposition}
\newtheorem{remark}{Remark}

\newcounter{unnumber}

\newcommand{\R}{\mathbb{R}}%
\DeclareMathOperator*\inte{int}%
\DeclareMathOperator*\cl{cl}%
\DeclareMathOperator*\epi{epi}%
\DeclareMathOperator*\dom{dom}%

\DeclareMathOperator*\B{\overline{\R}}%
\DeclareMathOperator*\pr{pr}%

\title{On a zero duality gap result in extended monotropic programming}

 \author{Radu Ioan Bo\c{t} \thanks
 {Faculty of Mathematics, Chemnitz University of Technology,
D-09107 Chemnitz, Germany, e-mail:
 radu.bot@mathematik.tu-chemnitz.de. Research partially supported by DFG (German Research Foundation), project WA 922/1-3.} \and Ern\"{o} Robert Csetnek
 \thanks {Faculty of Mathematics, Chemnitz University of Technology,
D-09107 Chemnitz, Germany, e-mail:
 robert.csetnek@mathematik.tu-chemnitz.de}}

\begin{document}
\maketitle

\noindent \textbf{Abstract.} In this note we correct and improve a zero duality gap result in extended monotropic
programming given by Bertsekas in \cite{bertsekas}.\\



\section{Preliminaries}

In this paper we deal with the \textit{extended monotropic programming problem} (for the origins of which we refer to \cite{rock-monotr-81, rock-monotr-98})
\begin{equation*}
\begin{array}{lrl}
(P) & \inf & \sum\limits_{i=1}^m f_i(x_i)\\
& \mbox{s.t.} & (x_1,...,x_m) \in S
\end{array}
\end{equation*}
and its dual problem
\begin{equation*}
\begin{array}{lrl}
(D) & \sup & \sum\limits_{i=1}^m -f_i^*(x_i^*),\\
& \mbox{s.t.} & (x_1^*,...,x_m^*) \in S^\bot
\end{array}
\end{equation*}
where $X_i$ are separated locally convex spaces, $f_i : X_i \rightarrow \overline{\R}$ are proper and convex functions, $i=1,...,m,$ and
$S \subseteq \prod_{i=1}^m X_i$ is a linear closed subspace such that $\prod_{i=1}^m \dom f_i \cap S \neq \emptyset$.

The same primal-dual pair has been recently investigated by Bertsekas in \cite{bertsekas} in the case $X_i = \R^{n_i}, n_i \geq 1$, $i=1,...,m$. In \cite[Proposition 4.1]{bertsekas}, under the supplementary assumption that the functions $f_i$ are lower semicontinuous on $\dom f_i$, a \textit{zero duality gap} result is stated for $(P)$ and $(D)$, provided that for every $(x_1,...,x_m) \in \prod_{i=1}^m \dom f_i \cap S$ and every $\varepsilon > 0$ the set
$$T(x,\varepsilon):= S^\bot+\prod_{i=1}^{m}\partial_{\varepsilon}f_i(x_i)$$
is closed. The proof of this statement, which represents the main result in that article, applies in an ingenious way the $\varepsilon$-descent method.

In this note we furnish first an example which shows that this zero duality gap statement is false and indicate the place where the error occurs. This will be the topic of the forthcoming section. In Section 3 we prove that under alternative, still weak, topological assumptions for the functions $f_i, i=1,...,m$, the zero duality gap statement in discussion turns out to be true and use to this aim some convex analysis specific techniques based on \textit{subdifferential calculus}, whereby a determinant role is played by a generalization of the \textit{Hiriart-Urruty--Phelps formula}. Recall that by zero duality gap we name the situation when $v(P) = v(D)$, where $v(P)$ and $v(D)$ denote the optimal objective values of the primal and dual problem, respectively.

In the following we introduce and recall some notions and results in order to make
the paper self-contained. Having a separated locally convex vector space $X$, we denote by $X^*$ its
topological dual space and assume throughout the paper that this is endowed with the weak$^*$ topology. By $\langle x^*, x\rangle=x^*(x)$ we denote the value of the continuous linear functional $x^*\in X^*$ at $x\in X$. Given a subset $U$ of $X$, by $\cl (U)$ we denote its \textit{closure}. By $\delta_U : X \rightarrow \overline\R = \R \cup \{\pm \infty\}$, defined by $\delta_U(x) = 0$ for $x \in U$ and $\delta_U(x) = +\infty$, otherwise, we denote its \textit{indicator function}, while by $\sigma_U : X^* \rightarrow  \overline\R$, defined by $\sigma_U(x^*) = \sup_{x \in U} \langle x^*,x \rangle$, its \textit{support function}. We call a set $K \subseteq X$ \textit{cone} if for all $\lambda \geq 0$ and all $k \in K$ one has $\lambda k \in K$. For a given cone $K \subseteq X$ we denote by $K^*=\{x^*\in X^*:\langle x^*, k\rangle \geq 0\ \forall k\in K\}$ its \textit{dual cone} and for $S \subseteq X$ a linear subspace we denote by $S^\bot = \{x^*\in X^*:\langle x^*, x\rangle = 0\ \forall x\in S\}$ its \textit{orthogonal space}. For $U,V \subseteq X$ two given sets, the projection operator $\pr_U : U \times V \rightarrow U$ is defined as $\pr_U(u,v) = u$ for all $(u,v) \in U \times V$.

Having a function $f:X\rightarrow \overline\R$ we use the
classical notations for its \textit{domain} $\dom f=\{x\in X:
f(x)<+\infty\}$, its \textit{epigraph} $\epi f=\{(x, r)\in X\times
\R: f(x)\leq r\}$ and its \textit{conjugate function}
$f^*:X^*\rightarrow \overline\R$, $f^*(x^*)=\sup\{\langle x^*,
x\rangle - f(x) : x\in X\}$. Regarding a function and its conjugate we have the \textit{Young-Fenchel inequality}
$f^*(x^*)+f(x)\geq \langle x^*, x\rangle$ for all $x\in X$ and
$x^*\in X^*$. We call $f$ \textit{proper} if
$f(x)>-\infty$ for all $x\in X$ and $\dom f \neq\emptyset$.

For ${\varepsilon}\geq 0$, if $f(x)\in \R$ the {\it
$\varepsilon$-subdifferential} of $f$ at $x$ is
$$\partial_{\varepsilon} f (x)=\{x^*\in X^*: f(y)-f(x)\geq \langle
x^*, y-x\rangle-{\varepsilon}\ \forall y\in X\},$$ while if
$f(x)=\pm \infty$ we take by convention $\partial_{\varepsilon} f
(x):=\emptyset$. We denote by $\partial f(x):=
\partial_{0} f(x)$ the \textit{(convex) subdifferential} of $f$ at $x$. The {$\varepsilon$-subdifferential} of $f$ at $x$
is always a convex and closed set.  If $f$ is a proper function, then for $x \in \dom f$, $x^* \in X^*$ and $\varepsilon\geq 0$
one has
$$f(x)+f^*(x^*)\leq \langle x^*, x\rangle + \varepsilon \Leftrightarrow x^*\in
\partial_{\varepsilon} f(x) \Rightarrow x \in \partial_{\varepsilon} f^*(x^*).$$
If $0\leq \varepsilon\leq \eta$ it holds
$\partial_{\varepsilon} f(x)\subseteq \partial_{\eta}f(x)$ and $\cap_{\mu > \varepsilon} \partial _{\mu}f(x) = \partial _{\varepsilon} f(x)$ for all $x\in X$. Assuming that $f$ is a proper and convex function and $x \in \dom f$, then (see, for instance, \cite[Theorem 2.4.4 (iii)]{Zal-carte}) $f$ is lower semicontinuous at $x$ if and only if $\partial_{\varepsilon} f(x) \neq \emptyset$ for all $\varepsilon > 0$. Therefore, if $f^*(x^*) \in \R$ and $\varepsilon > 0$ one has $\partial_\varepsilon f^*(x^*) \neq \emptyset$.

If $K$ is a nonempty cone, then $\delta_K^* = \sigma_K = \delta_{-K^*}$ and $\partial_{\varepsilon} \delta_K(0) = -K^*$ for all $\varepsilon \geq 0$, while,
if $S$ is a nonempty linear subspace, then $\delta_S^* = \sigma_S = \delta_{S^\bot}$ and $\partial_{\varepsilon} \delta_S(x) = S^\bot$ for all $\varepsilon \geq 0$ and all $x \in S$.

The \textit{lower semicontinuous hull} of $f:X\rightarrow \overline\R$ is the function $\cl f:X\rightarrow
\overline\R$ which has as epigraph $\cl(\epi f)$. One always has that $\dom f \subseteq \dom (\cl f) \subseteq \cl(\dom f)$ and $f^* = (\cl f)^*$. Assuming that $f$ is convex, $f^*$ is proper if and only if $\cl f$ is proper, the latter being a sufficient condition for $f^{**} = \cl f$. Given the proper functions $f,g :X\rightarrow \overline \R$,
their \textit{infimal convolution} is the function $f \square g :X\rightarrow
\overline\R$ , $(f \square g)(x)=\inf \{f(x-y) + g(y) : y \in X\}$. If $f,g :X\rightarrow \overline \R$ are proper, convex and lower semicontinuous functions with $\dom f \cap \dom g \neq \emptyset$, then one has the \textit{Moreau-Rockafellar formula} $(f+g)^* = \cl(f^*\square g^*)$ (see \cite{BW}). For the convex analysis notions and results introduced in this section we refer to \cite{EkTem, Zal-carte}.

We would like to close this section by pointing out that, for $g: \prod_{i=1}^m X_i \rightarrow \overline \R$, $g(x_1,...,x_m) = \sum_{i=1}^m f_i(x_i)$, the primal problem $(P)$ can be equivalently written as
$$\inf_{x=(x_1,...,x_m) \in \prod\limits_{i=1}^m X_i} [g(x) + \delta_S(x)] = -(g+\delta_S)^*(0).$$
Its \textit{Fenchel dual} problem is
$$\sup_{x^*=(x_1^*,...,x_m^*) \in \prod\limits_{i=1}^m X^*_i} [-g^*(x^*) - \delta_S^*(-x^*)] = -(g^* \square \delta_S^*)(0)$$
and, since for $x^*=(x_1^*,...,x_m^*) \in \prod_{i=1}^m X^*_i$ one has $g^*(x_1^*,^...,x_m^*) = \sum_{i=1}^m f_i^*(x_i^*)$,
this is further equivalent to
$$\sup_{(x_1^*,..,x_m^*) \in S^\bot} -\sum_{i=1}^m f_i^*(x_i^*),$$
being nothing else than the dual problem $(D)$. Thus one can notice that for the primal-dual pair in discussion we always have \textit{weak duality}, i.e. $v(P) \geq v(D)$.

\section{Examples}

In the beginning of this section we give the announced example, which shows that under the hypotheses considered in \cite{bertsekas} the duality statement \cite[Proposition 4.1]{bertsekas} may fail.

\begin{example}\label{berts-fail}\rm Consider the convex set $C = \{0\}\times[3,\infty) \cup\inte(\R^2_+)$ and define the functions $f_1 :\R^2\rightarrow\B$ by $f_1(u,v)= v + \delta_C(u, v)$ and $f_2 : \R \rightarrow \overline {\R}$, $f_2(w) = \delta_{\R_-}(w)$. We are in the case $m=2, n_1=2, n_2=1$.  We further take   $S=\{(u,v,w)\in\R^{3} :u=w\}$, which is a linear subspace of $\R^3$ and show that the assumptions of \cite[Proposition 4.1]{bertsekas} are fulfilled. The functions $f_1,f_2$ are proper and convex, $f_1$ is lower semicontinuous on $\dom f_1=C$, $f_2$ is lower semicontinuous (on $\R$) and the feasible set of the primal problem is $(\dom f_1\times\dom f_2)\cap S=\{0\}\times[3,\infty)\times\{0\}$.

Next we prove that for all $\varepsilon>0$ and all $a\geq 3$, the set $T((0,a,0),\varepsilon) = S^\bot+\partial_{\varepsilon}f_1(0,a)\times\partial_{\varepsilon}f_2(0)$ is closed. Let us fix some arbitrary elements $\varepsilon>0$ and $a\geq 3$. One can easily see that $S^\bot=\{(x^*,0,-x^*):x^*\in\R\}$ and $\partial_{\varepsilon}f_2(0)=\R_+$. We claim that \begin{equation}\label{eq-eps-sub-f-2}\partial_{\varepsilon}f_2(0,a)=\R_-\times \left[1-\frac{\varepsilon}{a},1\right].\end{equation} According to the definition of the $\varepsilon$-subdifferential, an element $(u^*,v^*)$ belongs to $\partial_{\varepsilon}f_2(0,a)$ if and only if \begin{equation}\label{ineq1}v-a\geq u^*u+v^*(v-a)-\varepsilon \ \forall (u,v)\in C.\end{equation} We show first that $\R_-\times \left [1-{\varepsilon}/a,1 \right ]\subseteq \partial_{\varepsilon}f_2(0,a)$. Take $u^*\leq 0$ and $v^*\in \left[1-{\varepsilon}/a,1 \right]$. Then for each $(u,v)\in C$ we get  $$u^*u+v^*(v-a)-\varepsilon\leq u^*u+v^*v-a+\varepsilon-\varepsilon=u^*u+(v^*-1)v+v-a\leq v-a,$$ hence $(u^*,v^*)\in\partial_{\varepsilon}f_2(0,a)$. For the opposite inclusion, take an arbitrary element $(u^*,v^*)\in \partial_{\varepsilon}f_2(0,a)$. One can easily derive from \eqref{ineq1} that \begin{equation}\label{ineq2}v-a\geq u^*u+v^*(v-a)-\varepsilon \ \forall (u,v)\in\R^2_+.\end{equation} From here one has that $u^*\leq 0$. By taking $u:=0$ in \eqref{ineq2} we obtain \begin{equation}\label{ineq3}(v^*-1)(v-a)\leq\varepsilon \ \forall v\geq 0,\end{equation} thus $v^*\leq 1$. For $v:=0$ in \eqref{ineq3} we get $a(v^*-1)\geq-\varepsilon$, that is $v^*\geq 1-{\varepsilon}/a$. In conclusion, \eqref{eq-eps-sub-f-2} holds. As a consequence we get $$T((0,a,0), \varepsilon)=\{(x^*,0,-x^*):x^*\in\R\}+\R_-\times \left [1-\frac{\varepsilon}{a},1 \right ]\times\R_+ = \R\times \left [1-\frac{\varepsilon}{a},1 \right]\times\R,$$ which is a closed set. Hence all the hypotheses of \cite[Proposition 4.1]{bertsekas} are fulfilled.

However, there is a nonzero duality gap between the primal-dual pair $(P)-(D)$. Indeed, $$v(P)=\inf_{(u,v,w)\in S}\{f_1(u,v)+f_2(w)\}=\inf_{(u,v,w)\in\{0\}\times[3,\infty)\times\{0\}} v=3,$$ while, $$v(D)=\sup_{(u^*,v^*,w^*)\in S^\bot}\{-f_1^*(u^*,v^*)-f_2^*(w^*)\}=\sup_{u^*\in\R}\{-f_1^*(u^*,0)-f_2^*(-u^*)\}$$$$=\sup_{u^* \leq 0}\Big\{-\sup_{(u,v)\in C}\{u^*u-v\}\Big\}=0.$$ Consequently, $v(D)<v(P)$, although the assumptions of \cite[Proposition 4.1]{bertsekas} are fulfilled.
\end{example}

Let us point out in the following where the error that occurred in \cite{bertsekas} comes from. The author claims that the formula $\sigma_{\partial_{\varepsilon}f(x)}=f_{\varepsilon}'(x,\cdot)$ is valid, where $f$ is a proper and convex function which is lower semicontinuous on $\dom f$, $x\in\dom f$ and $\varepsilon>0$ (cf. \cite[Section 3]{bertsekas}, see \cite[relation (15)]{bertsekas}). Here $f_{\varepsilon}'(x,y) = \inf_{\alpha > 0} (f(x+\alpha y) - f(x) + \varepsilon)/{\alpha}$ denotes the \textit{$\varepsilon$-directional derivative} of $f$ at $x$ in the direction $y \in X$. He decisively uses this formula in his argumentation, however, this formula holds in case $f$ is proper, convex and lower semicontinuous (on the whole space) (see \cite[Theorem 2.4.11]{Zal-carte} and \cite[p. 220]{Rock-carte}). Otherwise it can fail, as the following example shows.

\begin{example}\label{eps-dir-sigm-subdiff}\rm Consider $X$ a separated locally convex space and $K\subseteq X$ a nonempty convex cone which is not closed and define $f=\delta_K$. The function $f$ is proper, convex and lower semicontinuous on $\dom f=K$. Take $u\in\cl (K) \setminus K$ and $\varepsilon>0$. One can easily show that $f_{\varepsilon}'(0,u)=+\infty$ and $\partial_{\varepsilon}f(0)=-K^*$, hence $\sigma_{\partial_{\varepsilon}f(0)}(u)=\delta_{\cl (K)}(u)=0<f_{\varepsilon}'(0,u)$.
\end{example}

One of the main ingredients of the $\varepsilon$-descent method, on which the proof of the duality result \cite[Proposition 4.1]{bertsekas} relies, is \cite[Proposition 3.1]{bertsekas}. In its proof the formula discussed above is used, too. Let us recall this result: if $f_i:\R^n\rightarrow\B$ are proper and convex functions, $i=1,...,m$, and  $x\in\bigcap_{i=1}^{m}\dom f_i$ is a vector such that $f_i(x)=(\cl f_i)(x)$ for all $i=1,...,m$, then for all $\varepsilon>0$ the inclusion $$\partial_{\varepsilon}(f_1+...+f_m)(x)\subseteq \cl\big(\partial_{\varepsilon}f_1(x)+...+\partial_{\varepsilon}f_m(x)\big)$$ holds. We show in the following example that this is not always the case.

\begin{example}\label{prop3.1-bert}\rm Take $m=n=2$, $K=\inte(\R^2_+)\cup\{(0,0)\}$, $S=\R\times\{0\}$ and define the functions $f_1=\delta_K$ and $f_2=\delta_S$, which are proper and convex functions such that $\dom f_1\cap\dom f_2=\{(0,0)\}$. The vector $x=(0,0)$ satisfies the property $f_i(0,0)=(\cl f_i)(0,0)$, $i=1,2$. Take an arbitrary $\varepsilon>0$. One can show that $f_1+f_2=\delta_{\{(0,0)\}}$, hence $$\partial_{\varepsilon}(f_1+f_2)(0,0)=\R^2.$$ Further, $\partial_{\varepsilon}f_1(0,0)=-K^*=-\R^2_+$ and $\partial_{\varepsilon}f_2(0,0)=S^\bot=\{0\}\times\R$, thus  $$\cl\big(\partial_{\varepsilon}f_1(0,0)+\partial_{\varepsilon}f_2(0,0)\big)=\R_-\times\R.$$ Thus the assertion in \cite[Proposition 3.1]{bertsekas} does not hold in this particular case.
\end{example}

Finally, let us mention that the results stated in \cite{bertsekas} in finite dimensional spaces become valid if the functions $f_i$, $i=1,...,m$, are assumed to be proper, convex and lower semicontinuous on the whole space. In the next section we prove, by using a different technique than in \cite{bertsekas}, that these results remain true in a more general context and under weaker assumptions.

\section{Zero duality gap in extended monotropic programming}

For the beginning we provide a generalization of the Hiriart-Urruty--Phelps formula (see \cite[Theorem 2.1]{hiriart-ur-phelps} and \cite[Corollary 2.6.7]{Zal-carte}). We refer the reader to \cite[Theorem 13]{h-l-zal}, \cite[Proposition 2]{dinh-lopez-volle} and \cite[Theorem 4]{lopez-volle} for other generalizations of this result. The proof of the following theorem is an adaptation of the one given in \cite[Theorem 2.1]{hiriart-ur-phelps}.

\begin{theorem}\label{gen-hiriart-phelps} Let $X$ be a separated locally convex space and $f,g:X\rightarrow\B$ two convex functions such that $\cl f$ and $\cl g$ are proper and the following equality holds
\begin{equation}\label{cl-sum} \cl(f+g)=\cl f +\cl g.
\end{equation} Then for all $x\in X$ and all $\varepsilon\geq 0$ we have \begin{equation}\label{eps-subdiff-sum} \partial_{\varepsilon}(f+g)(x)=\bigcap_{\eta>0}\cl\left(\bigcup\limits_{\substack{\varepsilon_1,\varepsilon_2\geq
0\\
\varepsilon_1+\varepsilon_2=\varepsilon+\eta}}\Big(\partial_{\varepsilon_1}f(x)+\partial_{\varepsilon_2}g(x)\Big)\right).
\end{equation}
\end{theorem}

\noindent {\bf Proof.} Take $x\in X$ and $\varepsilon\geq 0$. The inclusion $"\supseteq"$ is always true (even in the case when \eqref{cl-sum} is not fulfilled), since $\bigcup_{\substack{\varepsilon_1,\varepsilon_2\geq
0\\
\varepsilon_1+\varepsilon_2=\varepsilon+\eta}}\Big(\partial_{\varepsilon_1}f(x)+\partial_{\varepsilon_2}g(x)\Big)\subseteq\partial_{\varepsilon+\eta}(f+g)(x)$. Take now an arbitrary element $x_0^*\in\partial_{\varepsilon}(f+g)(x)$. This is equivalent to \begin{equation}\label{def-eps} (f+g)^*(x_0^*)+(f+g)(x)\leq\langle x_0^*,x\rangle+\varepsilon. \end{equation} We apply the Moreau-Rockafellar formula to the proper, convex and lower semicontinuous functions $\cl f$ and $\cl g$ and obtain (by using \eqref{cl-sum}) \begin{equation}\label{mor-rock-cl}(f+g)^*=\big(\cl(f+g)\big)^*=(\cl f+\cl g)^*=\cl\big((\cl f)^*\Box(\cl g)^*\big)=\cl(f^*\Box g^*).\end{equation} Thus by \eqref{def-eps} and \eqref{mor-rock-cl} it holds $(\cl\phi)(x_0^*)\leq r$, where $\phi:X^*\rightarrow\B$ is defined by $\phi(x^*)=(f^*\Box g^*)(x^*)-\langle x^*,x\rangle$ and $r:=\varepsilon-(f+g)(x) \in \R$. Let us fix an arbitrary $\eta >0$. The condition $(\cl\phi)(x_0^*)\leq r$ implies that \begin{equation}\label{x0-phi} x_0^*\in\cl\Big(\{x^*\in X^*:\phi(x^*)\leq r+\eta/2\}\Big).
\end{equation} Let us show that for $x^*\in X^*$ we have \begin{equation}\label{incl} \{x^*\in X^*:\phi(x^*)\leq r+\eta/2\}\subseteq \bigcup\limits_{\substack{\varepsilon_1,\varepsilon_2\geq
0\\
\varepsilon_1+\varepsilon_2=\varepsilon+\eta}}\Big(\partial_{\varepsilon_1}f(x)+\partial_{\varepsilon_2}g(x)\Big).\end{equation} Indeed, if $x^* \in X^*$ satisfies $\phi(x^*)-r\leq \eta/2$, then \begin{equation}\label{ineq-phi}\inf\limits_{\substack{x_1^*,x_2^*\in X^*\\x_1^*+x_2^*=x^*}}\Big\{f^*(x_1^*)+f(x)-\langle x_1^*,x\rangle+g^*(x_2^*)+g(x)-\langle x_2^*,x\rangle\Big\}<\varepsilon+\eta,\end{equation} hence there exist $x_1^*,x_2^*\in X^*$, $x_1^*+x_2^*=x^*$, such that \begin{equation}\label{ineq-f-g}f^*(x_1^*)+f(x)-\langle x_1^*,x\rangle+g^*(x_2^*)+g(x)-\langle x_2^*,x\rangle<\varepsilon+\eta.\end{equation} We define $\varepsilon_1:=f^*(x_1^*)+f(x)-\langle x_1^*,x\rangle$ and $\varepsilon_2:=\varepsilon+\eta-\big(f^*(x_1^*)+f(x)-\langle x_1^*,x\rangle\big)$. By using the Young-Fenchel inequality and \eqref{ineq-f-g} we easily derive that $\varepsilon_1,\varepsilon_2\geq 0$, $\varepsilon_1+\varepsilon_2=\varepsilon+\eta$, $x_1^*\in\partial_{\varepsilon_1}f(x)$ and $x_2^*\in\partial_{\varepsilon_2}g(x)$, hence \eqref{incl} holds. Combining \eqref{x0-phi} and \eqref{incl} we get the desired conclusion. \hfill{$\Box$}

\begin{remark}\label{cond-cl-f+g}\rm (i) Let us notice that the condition \eqref{cl-sum} is automatically fulfilled if we assume that $f$ and $g$ are lower semicontinuous.

(ii) If $f$ (or $g$) is finite and continuous at $x_0\in\dom f\cap\dom g$, then \eqref{cl-sum} holds (cf. \cite[Lemma 15]{h-l-zal}).

(iii) Let us mention that the condition \eqref{cl-sum} was used also by other authors (see \cite{h-l-zal, li-fang-lopez}) in order to generalize duality results or subdifferential formulae for convex functions which are not necessarily lower semicontinuous (see also \cite{dinh-lopez-volle, lopez-volle} for some nonconvex versions of these results).
\end{remark}

The formula of the $\varepsilon$-subdifferential of the infimal convolution of two functions, given in the proposition below, will play a decisive role in the proof of the main result of this section.

\begin{proposition}\label{prop-subdiff-inf-conv} (cf. \cite[Corollary 2.6.6]{Zal-carte}) Let $X$ be a separated locally convex space and $f_1,f_2: X \rightarrow \overline \R$ two proper and convex functions for which \begin{equation}\label{cond-subdiff-inf} \exists x^*\in X^*,\exists\alpha\in\R,\forall x\in X,\forall i\in\{1,2\}: f_i(x)\geq\langle x^*,x\rangle+\alpha.\end{equation} If $(f_1\Box f_2)(x)\in\R$ and $\varepsilon\geq 0$, then \begin{equation}\label{subdiff-inf-conv}\partial_{\varepsilon}(f_1\Box f_2)(x)=\bigcap_{\eta>0}\bigcup\limits_{\substack{y\in X,\varepsilon_1,\varepsilon_2\geq 0\\ \varepsilon_1+\varepsilon_2=\varepsilon+\eta}}\Big(\partial_{\varepsilon_1}f_1(x-y)\cap\partial_{\varepsilon_2}f_2(y)\Big).\end{equation}
\end{proposition}

\begin{remark}\label{rem-sudiff-inf-conv}\rm One can easily show that condition \eqref{cond-subdiff-inf} in the above statement is nothing else than $\dom f_1^*\cap\dom f_2^*\neq\emptyset$.
\end{remark}

Next we present the main result of the paper, which is a zero duality gap theorem for extended monotropic programming problems in infinite dimensional spaces stated under weak topological assumptions.

\begin{theorem}\label{dual-monotr} Let $X_i$ be separated locally convex spaces, $f_i : X_i \rightarrow \overline{\R}$ proper and convex functions, $i=1,...,m,$ $S \subseteq \prod_{i=1}^m X_i$ a linear closed subspace such that $\prod_{i=1}^m \dom f_i \cap S \neq \emptyset$ and $g:\prod_{i=1}^{m} X_i\rightarrow\B$ defined by $g(x_1,...,x_m)=\sum_{i=1}^{m}f_i(x_i)$. Suppose further that $\cl f_i$, $i=1,...,m$, are proper functions and $g(x)=(\cl g)(x)$ for all $x\in\dom(\cl g)\cap S$. If for all $(x_1,...,x_m)\in \prod_{i=1}^{m} \dom f_i \cap S$ and all $\varepsilon>0$ the set $$S^\bot+\prod_{i=1}^{m}\partial_{\varepsilon}f_i(x_i)$$
is closed, then $v(P)=v(D)$.
\end{theorem}

\noindent {\bf Proof.} If $v(P)=-\infty$, then $v(P)=v(D)$ holds by weak duality, therefore we consider in the following the case $v(P)\in\R$ (that $v(P)<+\infty$ is guaranteed by the feasibility assumption). By the hypotheses one has that $(\cl g)(x_1,...,x_m)=\sum_{i=1}^{m}(\cl f_i)(x_i)$ for all $(x_1,...,x_m)\in \prod_{i=1}^{m} X_i$, thus $\cl g$ is a proper function. Let us show now that \begin{equation}\label{cl-sum-s-g}\cl(\delta_{S}+g)=\delta_{S}+\cl g.\end{equation} The inequality "$\geq$" is always fulfilled, hence it is enough to prove that $\cl(\delta_{S}+g)(x) \leq (\delta_{S}+\cl g)(x)$ for all $x\in\dom(\cl g)\cap S$. Taking an arbitrary $x\in\dom(\cl g)\cap S$ we have $$\cl(\delta_{S}+g)(x)\leq(\delta_{S}+g)(x)=g(x)=(\delta_{S}+\cl g)(x)\leq \cl(\delta_{S}+g)(x),$$ thus \eqref{cl-sum-s-g} holds. The following inclusions (which can be proved by using the Young-Fenchel inequality) will be useful in what follows \begin{equation}\label{incl-eps-2eps}\partial_{\varepsilon}g(x_1,...,x_m)\subseteq \prod_{i=1}^{m}\partial_{\varepsilon}f_i(x_i)\subseteq \partial_{2\varepsilon}g(x_1,...,x_m) \ \forall (x_1,...,x_m)\in \prod_{i=1}^{m} X_i \  \forall\varepsilon\geq 0.\end{equation} We prove next that $(\delta_{S}^*\Box g^*)(0)\in\R$ and $\partial_{\varepsilon}(\delta_{S}^*\Box g^*)(0)\neq\emptyset$ for all $\varepsilon>0$.

Take an arbitrary $\varepsilon>0$. Since $(\delta_{S}+g)^*(0)=-v(P)\in\R$, we get $\partial_{\varepsilon/2}(\delta_{S}+g)^*(0)\neq\emptyset$. Let us choose an arbitrary $\overline{x}\in \partial_{\varepsilon/2}(\delta_{S}+g)^*(0)$. Thus $$(\delta_{S}+g)^*(0)+(\delta_{S}+g)^{**}(\overline{x})\leq\varepsilon/2.$$ Since $\cl(\delta_{S}+g)$ is a proper function, we get $$(\delta_{S}+g)^*(0)+ \cl(\delta_{S}+g)(\overline{x})\leq\varepsilon/2,$$ which implies $$(\delta_{S}+g)^*(0)+\delta_{S}(\overline{x})+(\cl g)(\overline{x})\leq\varepsilon/2,$$ hence $\overline{x}\in\dom(\cl g)\cap S$. Consequently, $(\cl g)(\overline{x})=g(\overline{x})$, $\overline{x}\in\dom g\cap S$ and  $$(\delta_{S}+g)^*(0)+\delta_{S}(\overline{x})+ g(\overline{x})\leq\varepsilon/2,$$ which is nothing else than $0\in\partial_{\varepsilon/2}(\delta_{S}+g)(\overline{x})$. Take an arbitrary $\eta>0$. We further apply Theorem \ref{gen-hiriart-phelps} and obtain $$\partial_{\varepsilon/2}(\delta_{S}+g)(\overline{x})\subseteq \cl\left(\bigcup\limits_{\substack{\varepsilon_1,\varepsilon_2\geq
0\\
\varepsilon_1+\varepsilon_2=(\varepsilon+\eta)/2}}\Big(\partial_{\varepsilon_1}\delta_{S}(\overline{x})+\partial_{\varepsilon_2}g(\overline{x})\Big)\right).$$
Since for $\varepsilon_1\geq 0$ we have $\partial_{\varepsilon_1}\delta_{S}(\overline{x})= S^\bot$, we get $$\partial_{\varepsilon/2}(\delta_{S}+g)(\overline{x})\subseteq\cl\left(\bigcup\limits_{\substack{\varepsilon_2\geq
0\\
\varepsilon_2\leq(\varepsilon+\eta)/2}}\Big(S^\bot+\partial_{\varepsilon_2}g(\overline{x})\Big)\right)=\cl\big(S^\bot+\partial_{(\varepsilon+\eta)/2}g(\overline{x})\big).$$ If we consider $\overline{x}=(\overline{x}_1,...,\overline{x}_m)$, where $\overline{x}_i\in X_i$, $i=1,...,m$, by \eqref{incl-eps-2eps} we have $$\cl\big(S^\bot+\partial_{(\varepsilon+\eta)/2}g(\overline{x})\big)\subseteq \cl\Big(S^\bot+\prod_{i=1}^{m}\partial_{(\varepsilon+\eta)/2}f_i(\overline{x}_i)\Big)$$
$$=S^\bot+\prod_{i=1}^{m}\partial_{(\varepsilon+\eta)/2}f_i(\overline{x}_i)\subseteq S^\bot+\partial_{\varepsilon+\eta}g(\overline{x}),$$ where we used the fact that the set $S^\bot+\prod_{i=1}^{m}\partial_{(\varepsilon+\eta)/2}f_i(\overline{x}_i)$ is closed. All together it follows that $0\in S^\bot+\partial_{\varepsilon+\eta}g(\overline{x})$. Hence there exists $y_0^*\in \partial_{\varepsilon+\eta}g(\overline{x})$ such that $-y_0^*\in S^\bot$. Thus $-y_0^*\in\partial\delta_{S}(\overline{x})$ and $y_0^*\in \partial_{\varepsilon+\eta}g(\overline{x})$ and from here we deduce that $\overline{x}\in\partial(\delta_{S}^*)(-y_0^*)\cap\partial_{\varepsilon+\eta}g^*(y_0^*)$. Hence $0=-y_0^*+y_0^*\in\dom\delta_{S}^*+\dom g^*=\dom(\delta_{S}^*\Box g^*)$ and (since $\eta>0$ is arbitrary) $$\overline{x}\in \bigcap_{\eta>0}\bigcup\limits_{\substack{y^*,\varepsilon_1,\varepsilon_2\geq 0\\ \varepsilon_1+\varepsilon_2=\varepsilon+\eta}}\Big(\partial_{\varepsilon_1}\delta_{S}^*(-y^*)\cap\partial_{\varepsilon_2}g^*(y^*)\Big).$$ As $\dom(\cl g)\cap S\neq\emptyset$, the condition \eqref{cond-subdiff-inf} (applied for $f_1=\delta_{S}^*$ and $f_2=g^*$) is fulfilled (see also Remark \ref{rem-sudiff-inf-conv}). The situation $(\delta_{S}^*\Box g^*)(0)=-\infty$, which would imply that $(\delta_{S}^*\Box g^*)^*=\delta_{S}+\cl g$ is identically $+\infty$, is not possible. Therefore, $(\delta_{S}^*\Box g^*)(0)\in\R$ and by Proposition \ref{prop-subdiff-inf-conv} we get $\overline{x}\in \partial_{\varepsilon}(\delta_{S}^*\Box g^*)(0)$.

Hence $\partial_{\varepsilon}(\delta_{S}^*\Box g^*)(0)\neq\emptyset$ for all $\varepsilon>0$. As $\delta_{S}^*\Box g^*$ is a proper and convex function and $0\in\dom(\delta_{S}^*\Box g^*)$, this implies that $\delta_{S}^*\Box g^*$ is lower semicontinuous at $0$. As in  \eqref{mor-rock-cl} (relation \eqref{cl-sum-s-g} holds) it follows that $(\delta_{S}+g)^*(0)=(\delta_{S}^*\Box g^*)(0)$ or, equivalently, $v(P)=v(D)$ and the proof is complete.\hfill{$\Box$}

\begin{remark}\label{lsc}\rm (i) Let us notice that in case the functions $\cl f_i$, $i=1,...,m$, are proper, the condition $g(x)=(\cl g)(x)$ for all $x\in\dom(\cl g)\cap S$ is satisfied if we assume that for all $i=1,...,m$, $f_i(x_i)=(\cl f_i)(x_i)$ for all $x_i\in\dom(\cl f_i)\cap\pr_{X_i}S$.

(ii) If the functions $f_i$ are lower semicontinuous on $X_i$, $i=1,...,m$, then the topological assumptions in Theorem \ref{dual-monotr}, namely that $\cl f_i$ are proper for $i=1,...,m$, and $g(x)=(\cl g)(x)$ for all $x\in\dom(\cl g)\cap S$ are obviously fulfilled.

(iii) We refer to \cite[Section 4.1]{bertsekas} for conditions which guarantee that for all $(x_1,...,x_m)\in \prod_{i=1}^{m} \dom f_i \cap S$ and all $\varepsilon>0$ the set $S^\bot+\prod_{i=1}^{m}\partial_{\varepsilon}f_i(x_i)$ is closed.
\end{remark}

\end{document}